 \documentclass[final,3p,times,12pt]{elsarticle}

\usepackage{times,amsmath,amssymb,dsfont,graphicx,color}

\journal{Operations Research Letters}

\def\Xc{\mathcal{X}}

\def\Ex{\textsc{E}}
\def\Pr{\textsc{P}}
\def\real{\mathbb{R}}

\def\bit{\begin{itemize}}
\def\eit{\end{itemize}}
\newcommand{\argmax}{\mathop{\mathrm{argmax}}}
\newcommand{\maxm}{\mathop{\mathrm{maximize}}}
\newcommand{\sbjt}{\mathop{\mathrm{subject\ to}}}
\newcommand{\Trns}{\top}
\newtheorem{theorem}{Theorem}
\newtheorem{corollary}{Corollary}
\newproof{IEEEproof}{Proof}

\begin{document}
\sloppy

\begin{frontmatter}

\title{On Bellman's principle with inequality
constraints}

\author[csuece]{Edwin K. P. Chong\corref{cor1}}
\ead{edwin.chong@colostate.edu}
\author[numerica]{Scott A. Miller}
\ead{scott.miller@numerica.us}
\author[numerica]{Jason Adaska}
\ead{jason.adaska@numerica.us}

\address[csuece]{Colorado State University, 1373 Campus Delivery, Fort
Collins, CO 80523, USA}
\address[numerica]{Numerica Corporation, 4850 Hahns Peak Drive, Suite
200, Loveland, CO 80538, USA}

\cortext[cor1]{Corresponding author}



\begin{abstract}
We consider an example by Haviv (1996) of a constrained
Markov decision process that, in some
sense, violates Bellman's principle. We resolve this
issue by showing how to preserve a form of Bellman's
principle that accounts for a change of constraint at states that
are reachable from the initial state.

\end{abstract}

\end{frontmatter}

\section{Introduction}

The most celebrated result in Markov decision process (MDP) theory is 
\emph{Bellman's optimality principle}, which can be stated as 
follows. (We assume that the reader is already generally familiar
with MDPs.) Let $X_t$ be the state at (discrete) time $t$ and
$r(X_t,a)$ the reward received if action $a$ is taken at state
$X_t$ (the stagewise reward). 
Let $V^*(x)$ be optimal cumulative reward starting at state $x$.
Then, Bellman's principle states that for each time $t$,
\[
V^*(X_t) = \max_a\{r(X_t,a) + \Ex_{X_t,a}[V^*(X_{t+1})]\}
\]
where $X_{t+1}$ is the random next state with distribution depending on
$X_t$ and $a$. Moreover, replacing max by argmax on the right-hand side
gives the optimal action at $X_t$ (i.e., it characterizes the optimal
policy). But Bellman's principle is more than just an equation---it
embodies an idea that has become almost fundamentally axiomatic in
Markov decision theory. This idea is that the optimal policy solves the
optimization problem not just at the initial state $X_0=x$ but also at
all states reachable from it.

In this paper, we consider MDPs with explicit constraints.
Such constrained MDPs have been studied for at
least a couple of decades and continues to draw interest (see, e.g., 
\cite{RoV89}--\cite{ChF07}).
We are interested here in a particular paper by Haviv \cite{Hav96}, who
raises an issue that has not been addressed in the literature.
Basically, Haviv constructs an example of a constrained MDP in which
the optimal policy starting at the initial state $x$ is no longer
optimal at states other than $x$, not even at a state $y$ that is
reachable from $x$. He laments that this means that Bellman's principle
is violated.

We will explore Haviv's issue thoroughly. In particular, we will show
that there is some preservation of Bellman's principle, provided we
account for the fact that some of the
``slackness'' in the constraint is
spent in going from $x$ to a reachable $y$. So, if we consider the
optimal policy $\pi^*$ starting at state $x$, the optimality of 
$\pi^*$ at state $y$ is with respect to a \emph{different} problem, one
where the constraint is modified with the ``residual slackness.''
In analyzing Haviv's problem, we will present some known results,
some new results, and some related examples along the
way to help us understand and resolve the problem. Our analysis
highlights the important maxim that when imposing constraints on a
decision problem, the constraints should apply only to those things
over which we have control.

\section{Haviv's Problem}

\begin{figure}
\begin{center}
\includegraphics[width=3in]{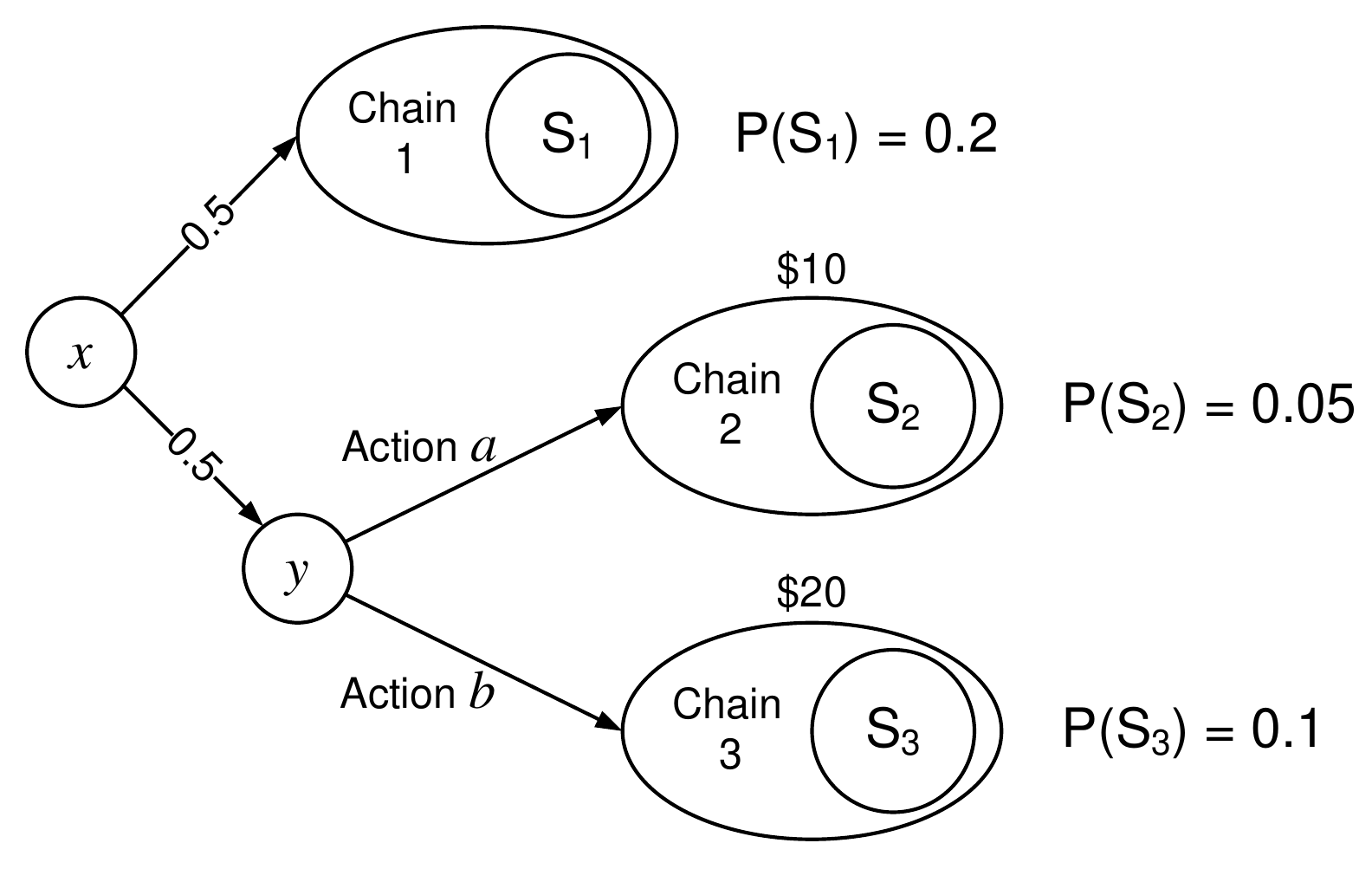}
\caption{Haviv's example}
\label{fig:havivf1}
\end{center}
\end{figure}

In \cite{Hav96}, Haviv gives an example (reproduced in 
Fig.~\ref{fig:havivf1}) in which he shows
that, given an optimal policy for an optimization problem starting
at some state $x$, the policy  is not optimal with respect to the
same problem starting at a reachable state $y$. The structure of the problem 
is a \emph{multichain} MDP with initial state $x$, which is transient. 
There are three recurrent 
subchains that could be reached from $x$. There is no reward for being in
chain~1, while the stagewise reward is $\$10$ at every state in
chain~2 and 
$\$20$ at every 
state in chain~3. The constraint is that the expected frequency of visits to 
states in $S=S_1\cup S_2\cup S_3$ must not exceed $0.125$ (think of
states in $S$ as the ``bad'' states). While in chain~$i$ 
($i=1,2,3$), the frequency of visits to $S_i$ is as shown in
Fig.~\ref{fig:havivf1} (e.g., $0.2$ for $S_1$). There is only one
state in which an action decision must be made: In state $y$, we can
choose either action $a$ or $b$.

A quick examination of Haviv's problem shows that there is
only one feasible policy: At state $y$, select action $a$. 
Selecting action $b$ at state $y$ would violate the
constraint, because the resulting Markov chain would visit states in
$S$ with frequency $0.5(0.2+0.1)=0.15$. However,
if the starting state were $y$, we would want to pick action $b$,
because this leads to chain~3 where the stagewise reward exceeds
that of chain~2, and the frequency of visits to $S$ in chain~3
($S_3$) is $0.1$, which does not exceed the constraint of $0.125$.
As noted before, this leads to Haviv's lament---Bellman's 
principle is violated, because the
optimal policy starting at state $x$ is no longer optimal
starting at state $y$, even though $y$ is reachable from $x$.
As Haviv points out in \cite{Hav96} and we will emphasize again 
later, the issue is related to the multichain nature of the
example: that there are transient states and recurrent subchains
that are not reachable from each other.

More specifically, Haviv's problem illustrates that as far as
optimality of a policy is concerned, arriving at state $y$ from $x$ 
is different from starting at $y$. From this point of view,
the issue raised by Haviv appears to be related to that of 
\emph{time consistency} in risk averse multistage stochastic
programming, identified in a recent paper by Shapiro \cite{Sha09}.
The same issue is also discussed in the economics literature on
multistage decision problems arising in dynamic portfolios; see, e.g.,
\cite{BaC10}, \cite{Kih09}. This issue has been recognized for some
time in the context of time-varying preferences \cite{Sto56},
\cite{Pol68} and game-theoretic formalisms of such changing
tastes \cite{PeY73}, \cite{Gol80}.

In lamenting the violation of Bellman's principle, Haviv quotes
Denardo \cite{Den82} on the principle of optimality:
``An optimal policy has the property that
whatever the initial node (state) and initial arc (decision) are,
the remaining
arcs (decisions) must constitute an optimal policy with regard to
the node (state) resulting from the first transition.''
But, we ask, must the policy be optimal with respect to the \emph{same}
problem? Indeed,
Denardo concedes that: ``The term \emph{principle} of
optimality is, however, somewhat misleading; it suggests that this
is a fundamental truth, not a consequence of more primitive
things.''

We will show that for a constrained MDP, the optimal policy starting at
one state is optimal with respect to a problem with a \emph{modified}
constraint at each reachable state.  Basically, in going from state $x$
to $y$, we have ``spent'' some of the constraint, so the ``residual''
constraint is reduced.  We submit that this is not an unreasonable
predicament, and still satisfies Denardo's version of Bellman's
principle. Moreover, the articulation of Bellman's principle we derive
here is a consequence of basic optimality conditions (see
Theorems~\ref{thm:opt_av} and \ref{thm:optreach}), which we argue are
instances of ``more primitive things'' referred to by Denardo.

\section{Bellman's Principle}

\subsection{Notation}

We first provide a framework for analyzing MDPs with inequality
constraints, of the kind that is considered by Haviv \cite{Hav96}. We have
to set this up more rigorously than the statement of Bellman's
equation in the last section, because: (1) we wish to incorporate
explicit inequality constraints; (2) we consider the case of
expected long-term average reward (where Bellman's equation looks
slightly different); and (3) we need
sufficient generality for multichain problems.  For this reason, we
need some formal notation:
\begin{itemize}
\item State space: $\Xc$, assumed countable.
\item State sequence: $\{X_t\} = \{X_0, X_1, X_2, \ldots\}$.
\item Stagewise reward: $r(x,a)\in\real$
\item Stagewise constraint: $c(x,a)\in\real^n$
\item If $x$ is a state and $a$ an action, we write
$\Ex_{x,a}$ for the conditional expectation given $(x,a)$. For example,
if $L^*:\Xc\to\real$ is a given function, then $\Ex_{X_0,a}[L^*(X_1)]$ 
means that $X_1$ is distributed according to the transition probability 
distribution given $(X_0,a)$, and $\Ex_{X_0,a}[L^*(X_1)]$ is the conditional 
expectation of $L^*(X_1)$ with respect to this distribution given 
$(X_0,a)$. If a policy $\pi$ is given, then instead of writing 
$\Ex_{X_0,\pi(X_0)}[L^*(X_1)]$, we simply write $\Ex_{X_0}^\pi[L^*(X_1)]$. 
Similarly, given a policy $\pi$ and an initial state $X_0$, 
the distribution of the Markov process $\{X_t\}$ is well defined, 
and we write $\Ex_{X_0}^\pi$ for the conditional expectation with respect to this distribution.
\item Similarly, for conditional probability given an initial state $X_0$ and policy $\pi$, we use the notation
$\Pr_{X_0,a}$ and $\Pr_{X_0}^\pi$. We use ``$\Pr_{X_0}^\pi$-a.s.'' to 
mean almost surely (with probability one) with respect to the probability
measure $\Pr_{X_0}^\pi$.
\item For a vector $x\in\real^n$, we write $x\geq 0$ to mean nonnegativity 
of each component.
\end{itemize}

\subsection{Optimal Policy}

Fix a state $x\in\Xc$ and set the initial state $X_0=x$.
Let
\[
V_T^\pi(x) = \frac{1}{T}\sum_{t=0}^{T-1} r(X_t,\pi(X_t))
\]
and
\[
W_T^\pi(x) = \frac{1}{T}\sum_{t=0}^{T-1} c(X_t,\pi(X_t)).
\]
The objective function is given by
\[
V^\pi(x) = \Ex_{X_0}^\pi\left[ \lim_{T\to\infty} V_T^\pi(x) \right],
\]
and the constraint function by
\begin{equation}
W^\pi(x) = \Ex_{X_0}^\pi\left[ \lim_{T\to\infty} W_T^\pi(x) \right].
\label{eqn:W}
\end{equation}
With this notation, the optimization problem given $X_0=x$ is as
follows:
\begin{align}
\maxm_\pi\ & V^\pi(x) \nonumber\\
\sbjt\ & W^\pi(x) \geq 0. \label{eqn:opt}
\end{align}
Note that this form of the problem is sufficiently general to cover
other inequality constraints: $W^\pi(x) \leq 0$, $W^\pi(x) \leq w(x)$, etc.

First, we give sufficient conditions under which a policy is optimal
with respect to (\ref{eqn:opt}).  Though stated formally, we provide
this result not to claim any novelty in it, but merely so that we can
use it as a rigorous platform on which to frame our analysis. Indeed,
similar results can be found in the book by Altman \cite{Alt98},
though not exactly in this form (which is constructed explicitly for
the convenience of our current purposes). We also provide a proof,
using only elementary and familiar arguments, similar to those in the
book by Ross \cite{Ross70}.

\begin{theorem}\label{thm:opt_av}
Fix $x\in\Xc$ and set the initial state $X_0=x$.
Suppose there exist a policy $\pi^*$, a vector $\mu\in\real^n$, a
constant $V^*(x)\in\real$, and
a bounded function $L^*:\Xc\to\real$ such that the following hold:
\begin{description}
\item[(A1)] $W^{\pi^*}(x) \geq 0$
\item[(A2)] $\mu\geq 0$
\item[(A3)] $\mu^\Trns W^{\pi^*}(x) = 0$
\item[(A4)] $V^*(x) + L^*(X_t) = \max_a \{r(X_t,a) + \mu^\Trns c(X_t,a) +
\Ex_{X_t,a}[L^*(X_{t+1})]\}$
for $t=0,1,\ldots\,$, $\Pr_{X_0}^{\pi^*}$-a.s.
\item[(A5)] $\pi^*(X_t) \in \argmax_a \{r(X_t,a) + \mu^\Trns c(X_t,a) +
\Ex_{X_t,a}[L^*(X_{t+1})]\}$ 
for $t=0,1,\ldots\,$, $\Pr_{X_0}^{\pi^*}$-a.s.
\end{description}
Then $\pi^*$ is optimal with respect to (\ref{eqn:opt}) and $V^{\pi^*}(x) = V^*(x)$.
\end{theorem}

\begin{IEEEproof}
Let $\pi$ be a feasible policy. (Note that $\pi^*$ is feasible by
assumption (A1).) Then by assumption (A4),
$\Pr_{X_0}^{\pi^*}$-a.s.\ for $t=0,1,\ldots$\,,
\begin{align*}
& V^*(x) + L^*(X_t) \\
&= \max_a \{r(X_t,a) + \mu^\Trns c(X_t,a) +
\Ex_{X_t,a}[L^*(X_{t+1})]\} \\
&\geq r(X_t,\pi(X_t)) + \mu^\Trns c(X_t,\pi(X_t)) + \Ex_{X_t}^\pi[L^*(X_{t+1 })]
\end{align*}
with equality if $\pi=\pi^*$
(by assumption (A5)). 
Now multiply throughout by $1/T$ and sum from $0$ to $T-1$
to obtain
\begin{align*}
&\frac{1}{T}\sum_{t=0}^{T-1}V^*(x) + L^*(X_t) \geq \\
&\frac{1}{T}\sum_{t=0}^{T-1} r(X_t,\pi(X_t)) 
+  \mu^\Trns c(X_t,\pi(X_t)) + \Ex_{X_t}^\pi[L^*(X_{t+1})],
\end{align*}
which can be written as
\begin{align*}
& V^*(x) + \frac{1}{T}L^*(X_0) \geq  \\
& V_T^\pi(x) + \mu^\Trns W_T^\pi(x)
+ \left(\frac{1}{T}\sum_{t=1}^{T-1}\Ex_{X_{t-1}}^\pi[L^*(X_t)]-L^*(X_t)\right)\\
&\mbox{\qquad} + \frac{1}{T}\Ex_{X_{T-1}}^\pi[L^*(X_T)].
\end{align*}
Next, take limits as $T\to\infty$,
take expectation $\Ex_{X_0}^\pi$,
use the boundedness assumption on $L^*$,
and use the fact that 
$\Ex_{X_0}^\pi[\Ex_{X_{t-1}}^\pi[L^*(X_t)]]=\Ex_{X_0}^\pi[L^*(X_t)]$ 
(for $t\geq 1$) to obtain
\[
V^*(x) \geq V^\pi(x) + \mu^\Trns W^\pi(x)
\]
with equality if $\pi=\pi^*$.
Because $\pi$ is feasible, $W^\pi(x)\geq 0$.
Hence, because $\mu\geq 0$ by assumption (A2),
\[
V^*(x) \geq V^\pi(x).
\]
Now, for $\pi=\pi^*$, we use assumption (A3) to obtain
\[
V^*(x) = V^{\pi^*}(x),
\]
and in particular $V^{\pi^*}(x)\geq V^\pi(x)$.
This completes the proof.
\hfill$\boxdot$
\end{IEEEproof}
\bigskip

The equation
\begin{align*}
& V^*(X_0) + L^*(X_t) = \\
& \max_a \{r(X_t,a) + \mu^\Trns c(X_t,a) + \Ex_{X_t,a}[L^*(X_{t+1})]\}
\end{align*}
has the resemblance of Bellman's equation for the unconstrained case. 
Indeed, we can think of this optimality condition for the
constrained case as Bellman's equation associated with an
unconstrained MDP with stagewise reward $r(X_t,a) + \mu^\Trns
c(X_t,a)$ (called the \emph{Lagrangian reward}). However, the
optimality conditions (A1--A5) include not only Bellman's equation, but 
also conditions (equation and inequalities) akin to
Karush-Kuhn-Tucker (KKT) conditions. This form of the optimality
condition benefits from the usual interpretation of the
\emph{multiplier} vector $\mu$ as a ``price'' vector, and suggests
the possibility of approaching the problem using duality principles
(though we do not pursue this line of approach any further
here).

\subsection{Optimality at Subsequent Reachable States}

We say that a state $y$ is \emph{reachable} from $x$ at time $t\in\{0,1,\ldots\}$ under policy $\pi$ if, given $X_0=x$, we have $\Pr_{X_0}^\pi\{X_t=y\}>0$.
We say that $y$ is \emph{reachable} from $x$ under $\pi$ if 
there exists $t\in\{0,1,\ldots\}$ such that it is reachable at $t$.

Next, we show that the sufficient conditions in Theorem~\ref{thm:opt_av}
are enough for the \emph{same} $L^*$ to satisfy the Bellman's equation at every 
reachable state.

\begin{theorem}
Fix $x\in\Xc$.
Suppose there exist a policy $\pi^*$, a vector $\mu\in\real^n$, a
constant $V^*(x)\in\real$, and a bounded function $L^*:\Xc\to\real$ such that 
assumptions (A1--A5) hold.
Then, for each state $y\in\Xc$ reachable from $x$ under $\pi^*$,
\begin{align*}
V^*(x) + L^*(y) 
&= \max_a \{r(y,a) + \mu^\Trns c(y,a) + \Ex_{y,a}[L^*(X')]\} \\
\pi^*(y)
&\in \argmax_a \{r(y,a) + \mu^\Trns c(y,a) + \Ex_{y,a}[L^*(X')]\},
\end{align*}
where $X'$ is distributed according to the transition distribution given $(y,a)$.
\end{theorem}

\begin{IEEEproof}
Suppose there is some state $y\in\Xc$ reachable from $X_0$ under $\pi^*$ 
such that $V^*(x) + L^*(y) \neq \max_a \{r(y,a) + \mu^\Trns c(y,a) +
\Ex_{y,a}[L^*(X')]\}$. Since $y$ is reachable, there is some 
$t\in\{0,1,2,\ldots\}$ such that $\Pr_{X_0}^{\pi^*}\{X_t=y\}>0$.
Let $A$ be the event that
$V^*(x) + L^*(X_t) \neq \max_a \{r(X_t,a) + \mu^\Trns c(X_t,a) +
\Ex_{X_t,a}[L^*(X_{t+1})]\}$. Then, by assumption, 
$\Pr_{X_0}^{\pi^*}(A) = 0$. However, we can also write
\begin{align*}
\Pr_{X_0}^{\pi^*}(A) 
&= \Pr_{X_0}^{\pi^*}(A|\{X_t=y\})\Pr_{X_0}^{\pi^*}\{X_t=y\} \\
& \mbox{\qquad} + \Pr_{X_0}^{\pi^*}(A|\{X_t\neq y\})\Pr_{X_0}^{\pi^*}\{X_t\neq y\}\\
&\geq \Pr_{X_0}^{\pi^*}(A|\{X_t=y\})\Pr_{X_0}^{\pi^*}\{X_t=y\} \\
&= \Pr_{X_0}^{\pi^*}\{X_t=y\}\\
&> 0,
\end{align*}
which is a contradiction. 

A similar argument yields
$\pi^*(y) \in \max_a \{r(y,a) + \mu^\Trns c(y,a) + \Ex_{y,a}[L^*(X')]\}$.
This completes the proof.
\hfill$\boxdot$
\end{IEEEproof}

The theorem above immediately implies that conditions (A4) and (A5)
(along sample paths) hold for any state $y$ reachable from
$x$.

\begin{corollary}\label{cor:bellman}
Fix $x\in\Xc$. 
Suppose there exist a policy $\pi^*$, a vector $\mu\in\real^n$, a
constant $V^*(x)\in\real$, and a bounded function $L^*:\Xc\to\real$ such that 
assumptions (A1--A5) hold.
Let $y\in\Xc$ be reachable from $x$ under $\pi^*$, and
suppose we set the initial state to be $X_0=y$.
Then, $\Pr_{X_0}^{\pi^*}$-a.s.\ for $t=0,1,2,\ldots\,$,
\begin{align*}
& V^*(x) + L^*(X_t)\\
 &= \max_a \{r(X_t,a) + \mu^\Trns c(X_t,a) +
\Ex_{X_t,a}[L^*(X_{t+1})]\} \\
& \pi^*(X_t)\\
 &\in \argmax_a \{r(X_t,a) + \mu^\Trns c(X_t,a) +
\Ex_{X_t,a}[L^*(X_{t+1})]\}.
\end{align*}
\end{corollary}
\bigskip

The result above shows that Bellman's \emph{equation} holds at all
reachable states. Specifically,
(A4) and (A5) hold for any state reachable from $x$ (with
the objective function value $V^*(x)$ and multiplier vector $\mu$).
But this is not enough to show that $\pi^*$ is \emph{optimal} at
any state reachable from $x$.
The key hurdle is feasibility (i.e., (A1)).
To be specific, suppose that state $y$ is reachable from $x$ under $\pi^*$. 
In general, it is not true that $\pi^*$ is optimal with respect to the problem
\begin{align*}
\maxm_\pi\ & V^\pi(y) \\
\sbjt\ & W^\pi(y) \geq 0.
\end{align*}
Indeed, it is easy to construct examples for which $\pi^*$ is not 
feasible for the above problem (e.g., Haviv's example \cite{Hav96}).
However, a modification to the constraint (which depends on $x$)
gives us an
optimization problem starting at $y$ for which $\pi^*$ is indeed
optimal, as stated below.

First, we need some additional notation. 
Given $X_0=x$, let $y$ be reachable from $x$ at time $t$ under $\pi^*$.
Define
\[
C_y(x) = -\Ex_{X_0}^{\pi^*}[W^{\pi^*}(X_t)|X_t\neq y]
\frac {\Pr_{X_0}^{\pi^*}\{X_t \neq y\}} {\Pr_{X_0}^{\pi^*}\{X_t =
y\}}.
\]
Note that if $\Ex_{X_0}^{\pi^*}[W^{\pi^*}(X_t)|X_t\neq y]\geq 0$, then
$C_y(x)\leq 0$. Moreover,
the smaller the value of $\Pr_{X_0}^{\pi^*}\{X_t = y\}$,
the larger the value of $|C_y(x)|$.

\begin{theorem}\label{thm:optreach}
Fix $x\in\Xc$.
Suppose there exist a policy $\pi^*$, a vector $\mu\in\real^n$, a
constant $V^*(x)\in\real$, and a bounded function $L^*:\Xc\to\real$ such that 
assumptions (A1--A5) hold with $X_0=x$.
Let $y$ be reachable from $x$ at time $t$ under $\pi^*$.
Then $\pi^*$ is optimal with respect to the problem
\begin{align*}
\maxm_\pi\ & V^\pi(y) \\
\sbjt\ & W^\pi(y) \geq C_y(x)
\end{align*}
and $V^{\pi^*}(y) = V^*(x) - \mu^\Trns C_y(x)$.
\end{theorem}

\begin{IEEEproof}
We have, given $X_0=x$,
\begin{align*}
& W^{\pi^*}(x) \\
&= \Ex_{X_0}^{\pi^*}\left[\lim_{T\to\infty} \frac{1}{T}\sum_{k=0}^{T-1} c(X_k,\pi^*(X_k))\right] \\
&= \Ex_{X_0}^{\pi^*}\left[\lim_{T\to\infty} \frac{1}{T}
\sum_{k=0}^{t} c(X_k,\pi^*(X_k)) \right.\\
& \left.\mbox{\qquad} + \frac{T-t}{T}
\left( \frac{1}{T-t} \sum_{k=t}^{T-1} c(X_k,\pi^*(X_k))\right)
\right] \\
&= \Ex_{X_0}^{\pi^*}\left[ \Ex_{X_{t}}^{\pi^*}\left[ 
\lim_{T\to\infty} \frac{1}{T-t} \sum_{k=t}^{T-1} c(X_k,\pi^*(X_k))
\right] \right] \\
&= \Ex_{X_0}^{\pi^*}\left[W^{\pi^*}(X_{t})\right] \\
&= W^{\pi^*}(y)\Pr_{X_0}^{\pi^*}\{X_t=y\}
+\Ex_{X_0}^{\pi^*}\left[W^{\pi^*}(X_t)|X_t\neq y\right]
\Pr_{X_0}^{\pi^*}\{X_t\neq y\} \\
& = \left(W^{\pi^*}(y) -C_y(x)\right)\Pr_{X_0}^{\pi^*}\{X_t=y\}.
\end{align*}
By assumption, $W^{\pi^*}(x)\geq 0$, which implies that
$W^{\pi^*}(y) - C_y(x)\geq 0$. Moreover, because $\mu^\Trns W^{\pi^*}(x)= 0$,
we have $\mu^\Trns (W^{\pi^*}(y) - C_y(x)) = 0$.

Now, set the initial condition $X_0=y$.
Define a new stagewise constraint function
$\bar{c}(\cdot,a) = c(\cdot, a) - C_y(x)$ (subtracting the same constant
for each $a$) and let $\bar{W}^{\pi^*}(y)=W^{\pi^*}(y)-C_y(x)$, which is
the expected average constraint function defined accordingly using $\bar{c}$, 
analogous to (\ref{eqn:W}).
From the above, we have  $\bar{W}^{\pi^*}(y)\geq 0$ and
$\mu^\Trns \bar{W}^{\pi^*}(y)= 0$.
By Corollary~\ref{cor:bellman}, $\Pr_{X_0}^{\pi^*}$-a.s.\ for $t=0,1,2,\ldots\,$,
\begin{align*}
& V^*(x) + L^*(X_t)\\
& = \max_a \{r(X_t,a) + \mu^\Trns c(X_t,a) +
\Ex_{X_t,a}[L^*(X_{t+1})]\}.
\end{align*}
Subtract $\mu^\Trns C_y(x)$ from both sides to obtain
\begin{align*}
& (V^*(x)- \mu^\Trns C_y(x)) + L^*(X_t) \\
& = \max_a \{r(X_t,a) + \mu^\Trns \bar{c}(X_t,a) + \Ex_{X_t,a}[L^*(X_{t+1})]\}.
\end{align*}
Finally, again by Corollary~\ref{cor:bellman}, $\Pr_{X_0}^{\pi^*}$-a.s.\ for $t=0,1,2,\ldots\,$,
\[
\pi^*(X_t) \in \argmax_a \{r(X_t,a) + \mu^\Trns c(X_t,a) +
\Ex_{X_t,a}[L^*(X_{t+1})]\}.
\]
Note that if we substitute $\bar{c}$ for $c$, the above still holds.
We can now apply Theorem~\ref{thm:opt_av} to obtain the desired result.
\hfill$\boxdot$
\end{IEEEproof}
\bigskip

The theorem above has the interpretation of Bellman's principle
for constrained problems. Recall that in the
unconstrained case, this principle states that if $\pi^*$ is an optimal 
policy for a problem starting at some state $x$, then it is also optimal for
a problem starting at any state $y$ reachable from $x$.  The main wrinkle
in the constrained case is that the constraint for the problem starting at
$y$ is different from starting at $x$, because we have to take into account
how the constraint function $W^{\pi^*}$ depends on 
other
states that can be
reached. Basically, $C_y(x)$ plays the role of a ``residual slackness'' of
the constraint at the reachable state $y$, after the sojourn from $x$ to $y$.

Note that if $C_y(x)>0$, then constraint is more stringent at $y$.
In this case, we
can interpret $C_y(x)$ as the constraint that is ``spent'' in
going from $x$ to $y$. On the other hand,
if $C_y(x)<0$, then we \emph{gain} some slackness in
going from $x$ to $y$ (i.e., constraint is \emph{less} stringent).

\subsection{Haviv's Example}

We can use Theorem~\ref{thm:optreach} to construct the optimization
problem starting at $y$ for which the given policy is indeed
optimal. We use the notation $1_S(\cdot)$ for the indicator function of $S$, 
so that $1_S(x)=1$ if $x\in S$, and $1_S(x)=0$ otherwise. We have:
\begin{itemize}
\item $c(\cdot,a) = 0.125 - 1_{S}(\cdot)$
\item $C_y(x) = -(0.125-0.2)(0.5/0.5) = 0.075$
\item $\bar{c}(\cdot,a) = 0.125 - 1_{S}(\cdot) - 0.075 = 0.05-1_{S}(\cdot)$
\end{itemize}
So, instead of needing the expected frequency of visits to
$S$ not to exceed $0.125$, at state $y$ the constraint
becomes $0.05$ (more stringent).
In other words, we ``spent'' $0.075$ of the constraint in going from 
$x$ to $y$, and
the ``residual'' constraint starting at state 
$y$ is that the frequency of visits to $S$ should not exceed $0.05$. 
In this case, clearly only action $a$ is feasible at $y$.

\section{Satisfying Haviv}

\subsection{Form of Constraint is Bad}

Would our modified form of Bellman's principle satisfy Haviv? We
suspect not. Haviv's point is that intuition dictates that the optimal
policy should pick action $b$ at state $y$, though he acknowledges that
such a policy would not be feasible with respect to the problem
(\ref{eqn:opt}). He therefore goes on to argue that the \emph{form} of
the constraint in (\ref{eqn:opt}) is problematic. 
The version of Bellman's principle in Theorem~\ref{thm:optreach} 
is not entirely satisfactory because, one could argue, the constraint
should not change depending on what happened in the past. 

This is related to the issue of \emph{time consistency} in
\cite{Sha09}. Shapiro \cite{Sha09} defines time consistency as ``the
requirement that at every state of the system our `optimal' decisions
should not depend on scenarios which we already know cannot happen in
the future.'' In Haviv's example, once we are in state $y$, we know
that we will not enter chain~1. Yet, it is the frequency of visits to
states in $S_1$ within chain~1 that causes action $b$ to be infeasible
at state $y$.

This seems to be a legitimate concern.
We further illustrate this concern below by applying our
result to a different example. This example is in contrast to
Haviv's, because it turns out that at reachable states that are
unlikely to be visited, the constraint might be unreasonably relaxed.

\subsubsection*{Example: Squander or save}

\begin{figure}
\begin{center}
\includegraphics[width=3in]{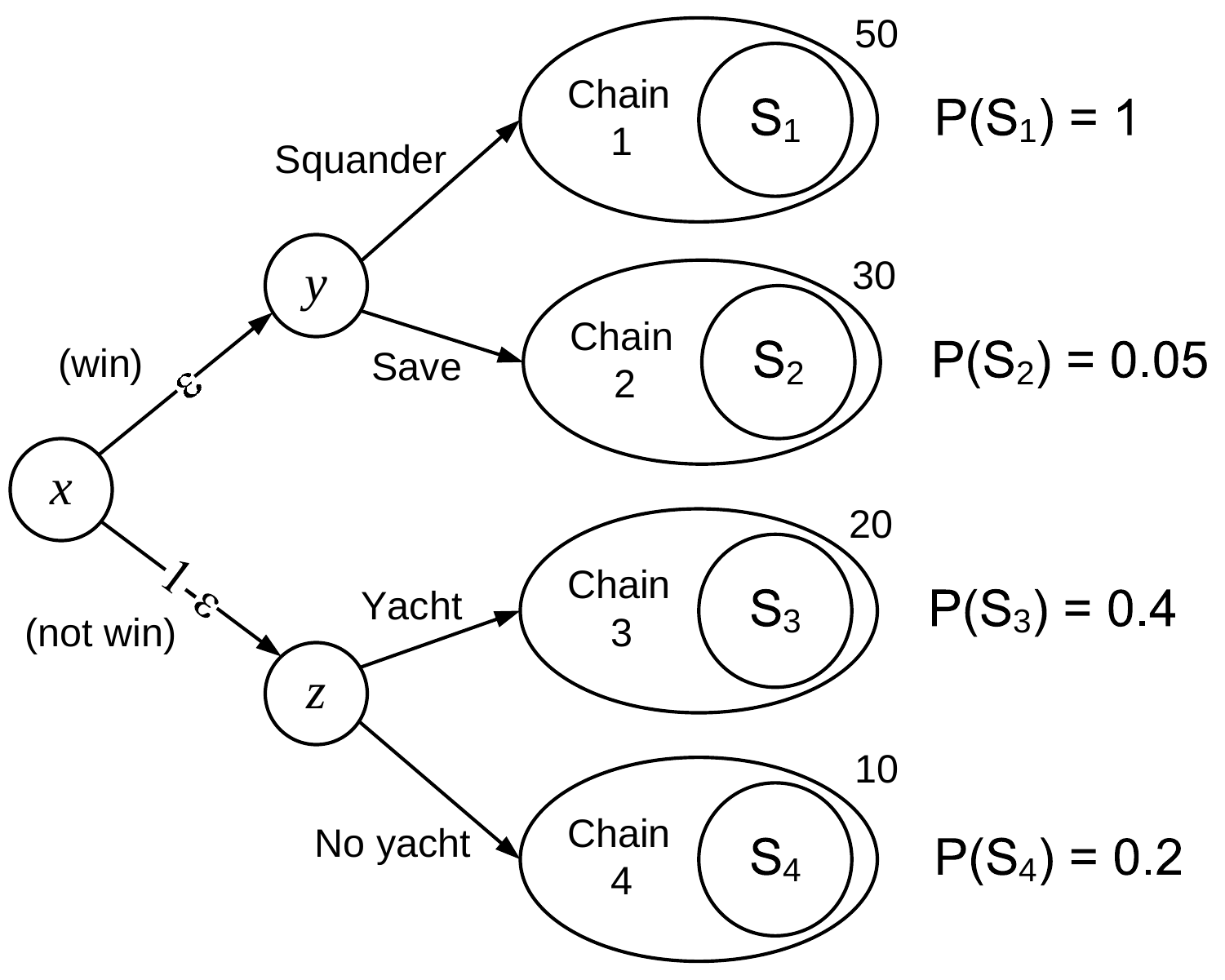}
\caption{Squander or save}
\label{fig:squander}
\end{center}
\end{figure}

Consider the problem in Fig.~\ref{fig:squander}. Starting at state
$x$, the process will go to either state $y$ or 
$z$
depending on
whether or not we win the lottery, respectively. The probability of
winning the lottery is (realistically) a small number $\varepsilon$,
as shown in the figure. If we do not win the lottery, we have the
choice of whether or not to buy a yacht. Depending on our choice, we
will end up in one of two possible subchains. In the unlikely event
that we \emph{do} win the lottery, we have the choice of whether or
not to squander all our money. Again, depending on this choice, we
end up in one of two possible subchains. Within each subchain, the
stagewise reward at all states is fixed at the value shown in the
figure (e.g., $50$ in chain~1). These reward values are meant to
signify the level enjoyment of life within these subchains.

In this example, the constraint is that the expected frequency of
visits to states in $S=S_1\cup S_2\cup S_3\cup S_4$ should not exceed
$0.3$. This constraint reflects the desire that we limit the
probability that we will go broke (have no money) before retiring. The
states in the problem that represent being broke are those in
$S=S_1\cup S_2\cup S_3\cup S_4$.  The frequency of visits to the
``bad'' states in each of the subchains is shown as $P(S_i)$,
$i=1,2,3,4$, in Fig.~\ref{fig:squander}.

It is clear that because $\varepsilon$ is taken to be very small, it is
overwhelmingly likely that we will enter state $z$, in which case we
cannot afford to buy a yacht---doing so would send us into chain~3,
where the frequency of visiting ``bad'' states is $0.4$, exceeding
$0.3$. But what about in state $y$, which corresponds to winning the
lottery?

In this problem, it turns out that $C_y(x) = 0.1(1-1/\varepsilon) <
0$.  So, depending on how small $\varepsilon$ is, $C_y(x)$ can be made
arbitrarily negative.  Specifically, for $\varepsilon<1/11$, the
optimal action at state $y$ is to squander. To be sure, it is not that
we can spend more if we win the lottery, but that because it is so
unlikely that we win, once we win we can do whatever we like without
violating the constraint.  This clearly illustrates that the form of
the constraint is problematic, as Haviv points out.

\subsection{Sample-Path Constraints}

How then can we resolve Haviv's problem?
Haviv advocates the use of \emph{sample-path} 
constraints, where we remove the expectation in the constraint
function in (\ref{eqn:W}) and
require instead that the inequality be satisfied with probability one. 
In our notation, this would correspond to, given $X_0=x$,
\[
 \lim_{T\to\infty} W_T^\pi(x) \geq 0\quad \mbox{$\Pr_{X_0}^\pi$-a.s.}
\]
It is clear that with such a constraint, a policy $\pi$ is feasible at 
$x$ if and only if feasible at each state reachable from $x$.

Note that this modification to the constraint immediately alleviates
the problem illustrated in Fig.~\ref{fig:squander}. In
contrast to the previous form of the constraint, it would
no longer be feasible to squander our money even if we win the
lottery.

\subsubsection*{Example: Yacht or not}

\begin{figure}
\begin{center}
\includegraphics[width=3in]{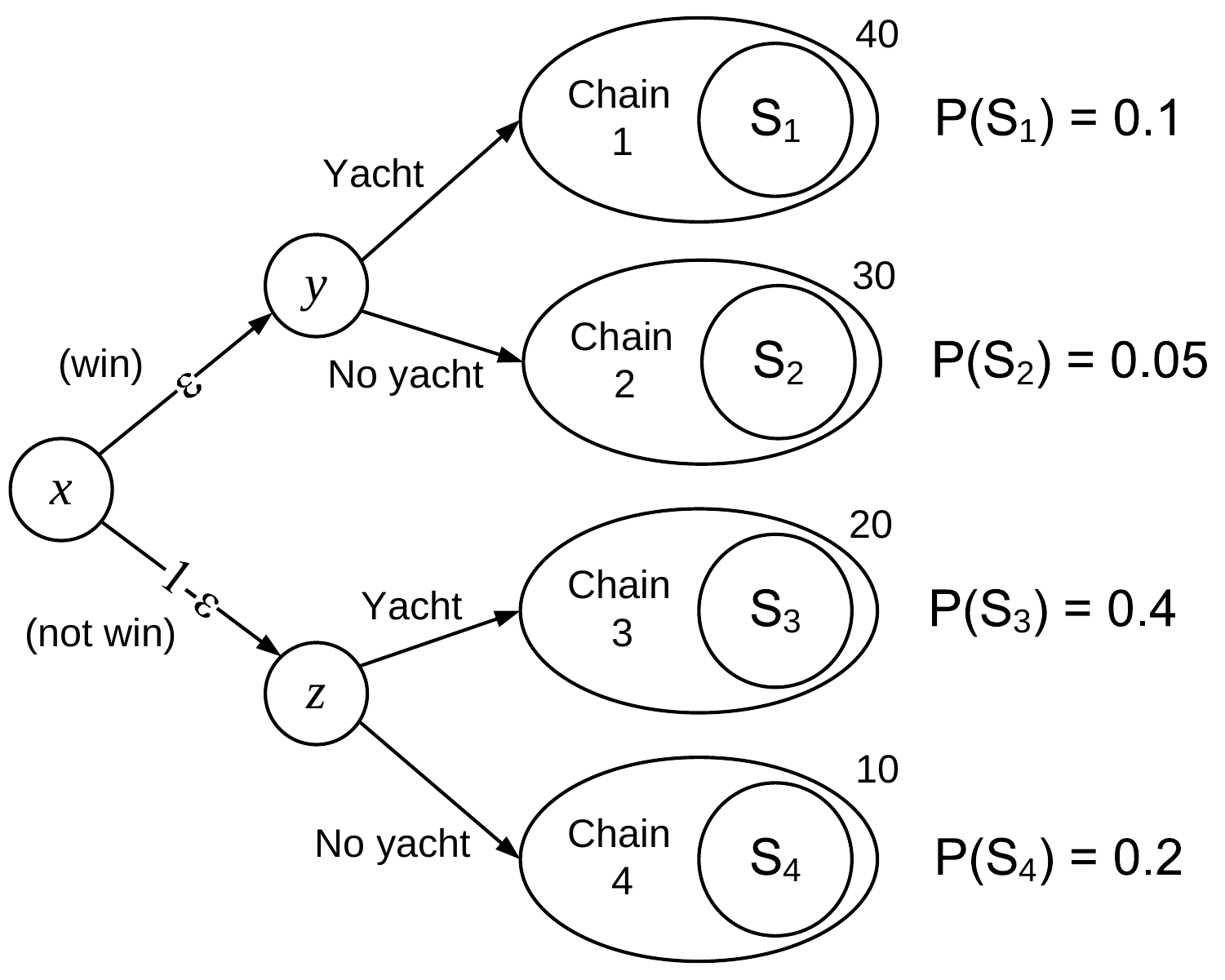}
\caption{Yacht or not}
\label{fig:lottery}
\end{center}
\end{figure}

To illustrate this point further, consider the problem in
Fig.~\ref{fig:lottery}, which is very similar to
Fig.~\ref{fig:squander} but simpler. In the current problem, again we
have the (unlikely) event of winning the lottery. However,
regardless of winning,
we have the decision of whether or not to buy a yacht.
Depending on whether or not we win the lottery and what decision we
make about the yacht, we will enter one of four subchains, wherein
there is some probability of going broke before retiring (as before,
these are shown as $P(S_i)$, $i=1,2,3,4$, in Fig.~\ref{fig:lottery}).
The stagewise reward values shown in the figure are again meant to
signify our enjoyment of life within these subchains.

As in the problem of Fig.~\ref{fig:squander}, if we impose sample-path
constraints, we will quickly arrive at the conclusion that in state
$z$, we cannot decide to buy a yacht because doing so would violate the
constraint in chain~3. However, in state $y$, where we have won the
lottery, we can in fact buy a yacht; doing so would not violate the
constraint. The optimal choice at state $y$ is indeed to buy a yacht,
leading to maximal enjoyment of life (within this example).

\subsection{Satisfying Haviv}

Suppose we make the modification from expected constraint to
sample-path constraint in Haviv's problem. Specifically, we now require
that, with probability one, the frequency of visits to $S$ not exceed
$0.125$. Then, \emph{no policy is feasible}, because there is a $0.5$
probability that the process will enter chain~1, where the frequency of
visits to $S$ is $0.2$. Haviv \cite{Hav96} does point this out, but does not
provide a resolution to it.

In other words, using a sample-path constraint does not resolve Haviv's
problem, because no policy would be feasible.  Moreover, for some class
of constrained MDPs (including Haviv's), sample-path constraints can be
converted to \emph{equivalent} expected constraints. These are what we
might call \emph{trans-policy decomposable} MDPs (see \cite{RoV89}).
Basically, to convert sample-path constraints into expected
constraints in such MDPs, we impose an expected constraint at each
subchain. For example, for subchain $C_3$, use the constraint function
$c(X_t,a)1_{\{X_t\in C_3\}}$ in the expected form of the constraint.

In the problems of Fig.~\ref{fig:squander} and Fig.~\ref{fig:lottery},
for example, we do not need to impose sample-path constraints; instead,
we can impose the usual (expected) form of the constraint in each of
the four subchains. If we do so, these problems would no longer suffer
from Haviv's problem, and the optimal policy would be equally optimal
at all reachable states without having to change the constraint.

The conversion of sample-path constraints into equivalent expected
constraints highlights an issue that is, at heart, what gives rise to
Haviv's problem: At state $y$
in Fig.~\ref{fig:havivf1}, we
have no control over whether the process will enter chain~1. Indeed,
this tells us that when imposing expected constraints in the subchains,
we should not impose them at \emph{all} subchains. In particular, we
should not impose any constraint at chain~1, because we have no control
(at $y$)
over whether or not we enter it. The constraints should be imposed only
at chains 2 and 3, which depend on a decision over which we have
control. If we do this, then even Haviv's original problem in
Fig.~\ref{fig:havivf1} would be resolved: The optimal policy with
respect to initial state $x$ is equally optimal (and feasible) at
state $y$, and would select action $b$ as desired.

Another way to express this observation is that constraints should be
imposed only on the \emph{consequence} of decisions, expressing conditions
on the desired (or undesired) impact of decisions once they are made.
In Haviv's example, expressing the constraint at state $x$ does not
properly reflect the impact of actions at $y$, which do not control
whether or not the system enters chain~1. The same would be true even
if we modify the example to include action choices at state
$x$ (e.g., we can control the probability of entering chain~1). The
constraint at state $x$ would still not properly reflect the impact of
actions at $y$, which do not control entry into chain~1, giving rise to
Haviv's lament.

\subsection*{Acknowledgement}

This work was supported in part by 
NSF Grant ECCS-0700559.

DISTRIBUTION STATEMENT A: Approved for public release; distribution is 
unlimited. Approved for public release 11-MDA-6364 (6 September 11).

\end{document}